\documentclass[11pt]{article}

\usepackage{amscd,amsmath, amsfonts, amssymb, fancyhdr, color}
\usepackage[a4paper]{geometry}

\newcommand{\ga}{\gamma}

\newcommand{\RR}{\mathbb{R}}

\numberwithin{equation}{section}

\newcommand{\ra}{{\:\longrightarrow\:}}
\newcommand{\Z}{{\mathbb Z}}
\newcommand{\C}{{\mathbb C}}
\newcommand{\R}{{\mathbb R}}

\newcommand{\cae}{{\cal E}}
 
\renewcommand{\bar}{\overline}
\renewcommand{\phi}{\varphi}
\renewcommand{\epsilon}{\varepsilon}

\newcommand{\Vol}{\operatorname{Vol}}

\newcommand{\Aut}{\operatorname{Aut}}

\renewcommand{\Im}{\operatorname{Im}}

\newcounter{Mycounter}[section]
\newcounter{lemma}[section]
\newcounter{claim}[section]
\newcounter{sublemma}[section]
\newcounter{corollary}[section]
\newcounter{theorem}[section]
\newcounter{conjecture}[section]
\newcounter{proposition}[section]
\newcounter{definition}[section]
\newcounter{example}[section]
\newcounter{remark}[section]
\newcounter{problem}[section]
\newcounter{question}[section]
\makeatletter
 
\@addtoreset{equation}{section}
\@addtoreset{footnote}{section}
\makeatother
 
\def\blacksquare{\hbox{\vrule width 5pt height 5pt depth 0pt}}
\def\endproof{\hfill\blacksquare}

 \newcommand{\samethanks}[1][\value{footnote}]{\footnotemark[#1]}
 
\begin{document}

\begin{center}
{\LARGE\bf
Holomorphic submersions of\\[3mm] locally conformally K\"ahler manifolds}\\[4mm]
{\large 
Liviu Ornea\footnote{Partially supported by CNCS UEFISCDI, project
number PN-II-ID-PCE-2011-3-0118.}, 
Maurizio Parton, and 
Victor Vuletescu\samethanks[1]}\\[4mm]

{\em Dedicated to the memory of Stere Ianu\c{s}}\\[4mm]
\end{center}




\begin{abstract}
A locally conformally K\"ahler (LCK) manifold is a complex manifold  covered by a K\"ahler manifold, with the covering group acting by homotheties. We show that if such a compact manifold $X$ admits a holomorphic submersion with positive dimensional fibers at least one of which is of K\"ahler type, then $X$ is globally conformally K\"ahler or biholomorphic, up to finite covers, to a Vaisman manifold ({\em i.e.} a mapping torus over a circle, with Sasakian fibre). As a consequence, we show that the product between a compact non-K\"ahler LCK and a compact K\"ahler manifold cannot carry a LCK metric.
\end{abstract}


\section{Introduction and statement of results}

{\em Locally conformally K\"ahler} (LCK) manifolds are Hermitian manifolds $(M,g,J)$ such that the fundamental two-form $\omega=g\circ J$ satisfies the integrability condition
$$d\omega=\theta\wedge\omega,\quad \text{for a {\em closed} one-form $\theta$},$$
where $\theta$ is called the Lee form.

This definition is known to be equivalent with a covering space $\tilde M$ of $(M,J)$ to carry a global K\"ahler metric $\Omega$ with respect to which the covering group $\Gamma$ acts by holomorphic homotheties (see \cite{GOPRVS}, and \cite{do}, \cite{ov11} for a recent survey). As such, the LCK structure defines a character associating to each covering transformation its scale factor:
\begin{equation}\label{ch}
\chi:\Gamma\ra\RR^+,\quad \chi(\ga)=\frac{\ga^*\Omega}{\Omega}.
\end{equation}

If $\theta$ is exact,  the manifold is called {\em globally conformally K\"ahler} (GCK) and is of K\"ahler type. 

In a LCK manifold, if $\theta$ is moreover  parallel with respect to the Levi-Civita connection of the LCK metric, the manifold is called Vaisman. Compact Vaisman manifolds are mapping tori over the circle with fibers isometric with a Sasakian manifold (see \cite{ov03}). The topology of compact Vaisman manifolds is very different from the topology of K\"ahler manifolds, {\em e.g.} their first Betti number is always odd.

Almost all  compact complex surfaces in class VII are LCK  and many of them ({\em e.g.} diagonal Hopf, Kodaira) are Vaisman (see \cite{be}, \cite{br}). In higher dimensions, main examples are diagonal Hopf manifolds (which are Vaisman), non-diagonal Hopf manifolds (non-Vaisman, see \cite{ov10}), Oeljeklaus-Toma manifolds (see \cite{ot, pv}).

On a Vaisman manifold $X$, the  Lee field $\theta^\sharp$ (the $g$-dual of $\theta$) is analytic and Killing and hence generates a complex, totally geodesic foliation $\mathcal{F}=\{\theta^\sharp, J\theta^\sharp\}$. 
If $\mathcal{F}$ is regular (and in this case the manifold $X$ itself is called \emph{regular}), then $X$ admits a holomorphic submersion (which is moreover a principal bundle map) over a projective orbifold. But, in general, very little is known about the existence of holomorphic submersions from compact LCK manifolds (papers like  \cite{iov}, \cite{mar} assume the existence of the submersion  and are mainly concerned by  the structure it imposes on the total space or on the base and by the geometry of the fibres).

In this note, we partially solve the existence problem. Our principal result is the following:

\hfill

\noindent{\bf Theorem.} {\em 
Let $X$ be a compact complex manifold which admits a homolorphic submersion  $\pi:X\ra B$ with positive-dimensional fibers. Assume one of the fibers of $\pi$ is of K\"ahler type. If $X$ has an LCK metric $g$, then $g$ is GCK or $X$ is biholomorphic to a  finite quotient of a Vaisman manifold.}

\hfill

This result is very general, it does not assume that the submersion relate in any way the Riemannian geometries of the total and base spaces (which is not even supposed to exist).

There is no natural  product construction in the category of LCK manifolds, because $\mathrm{CO}(m)\times \mathrm{CO}(n)\not \subset \mathrm{CO}(m+n)$. The following by-product of the Theorem (already proven differently in \cite{ts}) is therefore an useful information:

\hfill

\noindent{\bf Corollary 1.} {\em 
Let $X_1, X_2$ be compact regular Vaisman manifolds. Then $X_1\times X_2$ carries no LCK metric.
}

\hfill 

But more can be said. Applying the above Theorem to the projection of the first factor of a product $X\times Y$ where $X$ is a compact LCK (non-K\"ahler) manifolds and $Y$ is compact K\"ahler, one obtains:

\hfill

\noindent{\bf Corollary 2.} {\em The product of a compact LCK non-K\"ahler manifold with a compact K\"ahler manifold admits no  LCK metric.}

\section{Proof of the Theorem}

The main ingredient is the following ``lemma on fibrations''.

\hfill

\noindent{\bf Lemma.}\label{lemafib}
{\em Let $X$ be a compact complex manifold which admits a homolorphic submersion $\pi:X\ra B$ with positive-dimensional fibers. If $X$ has an LCK metric $g$ whose Lee form $\theta$ is (cohomologically) a pull-back, $[\theta]=\pi^*([\eta]), [\eta]\in H^1(B)$, then $g$ is GCK.}

\noindent {\bf Proof.} 
 The proof is basically the same as in \cite{ovv}, but since the statement is a little bit different, we include the details here.

First, let us fix some notations. If $M$ is any manifold and $\alpha\in H^1(M)$ is arbitrary, we will denote by $\alpha_*:H_1(M, \Z)\ra \R$ the morphism given by
$$\alpha_*([\gamma])=\int_\gamma \alpha.$$
Notice that in our setup we have 
$$\eta_*\circ \pi_*=\theta_*$$
where $\pi_*:H_1(X, \Z)\ra H_1(B, \Z)$ is the map induced at homology by $\pi.$

Moreover, we will denote by $M^{ab}$ its maximal abelian cover, whose fundamental group is just 
$[\pi_1(M), \pi_1(M)]$. Observe that the deck group of $M^{ab}$ over $M$ is $H_1(M, \Z).$

Now let $K=\ker (\eta_*);$ it is a subgroup of $H_1(B, \Z)$ 
so letting $\overline{B}=B^{ab}/K$ we see $H_1(\overline{B}, \Z)\cong H_1(B^{ab}, \Z)/K.$ In particular, the pull-back $\overline{\eta}$ of $\eta$ to  $\overline{B}$ is exact, since $\overline{\eta}_*\equiv 0.$ 

Now let $\overline{X}=\overline{B}\times_B X,$  {\em i.e.} 

$$
\begin{CD}
\bar X
@>>> X\\
@V\bar\pi VV @VV\pi V\\
\bar B
@>>>B
\end{CD}
$$
Then  $\overline{X}$ is a cover of $X$ and the fibers of the induced map $\overline{\pi}:\overline{X}\ra \overline{B}$ are the same as the fibers of $\pi$, thus $\overline{\pi}$ is proper as well. Let $\overline{\theta}$ be the pull-back of $\theta$ to $\overline{X}.$ Since $\overline{[\theta]}=\overline{\pi}^*([\overline{\eta}])$ we see that  $\overline{\theta}$ is also exact, as $\overline{\eta}$ is exact. This implies that the pull-back $\overline{g}$ of $g$ to $\overline{X}$ is globally conformal to a K\"ahler metric $\omega$.  

Assume now that $g$ is not GCK. Then there exists a deck transformation $\overline{\varphi}$ of $\overline{X}$ acting on $\omega$ by a non-isometric homothety:
\begin{equation}\label{rho}
\overline{\varphi}^*(\omega)=\varrho \cdot \omega, \qquad \rho\neq 1.
\end{equation}
Let $F$ be any fiber of $\overline{\pi}.$ Since $F$ is also a fiber of $\pi$, it is compact, hence its volume $\Vol_\omega(F)$ is finite.
Let $F'=\varphi(F).$ Then $F'$ is also a fiber  of $\overline{\pi}$ and since $\omega$ is K\"ahler we have
$$\Vol_\omega(F)=\Vol_\omega(F').$$
But from \eqref{rho} we get
$$\Vol_\omega(F)=\varrho^{\dim_\C(F)}\Vol_\omega(F'),$$
a contradiction. \hfill{$\blacksquare$}

\hfill


\noindent{\bf Proof of the Theorem.} 
We shall prove the following facts:
\begin{enumerate}
\item[{\bf 1)}]If the fibers are at least $2$-dimensional, then $g$ is GCK. 
\item[{\bf 2)}] If the fibers are $1$-dimensional and their genus is not $1$, $g$ is again GCK. And finally
\item[{\bf 3)}]  If the fibers are $1$-dimensional and their genus is  $1$, then $X$ is biholomorphic to a GCK manifold or to a finite quotient of a Vaisman manifold.
\end{enumerate}

To prove {\bf 1)}, let $F_0$ be a fiber of K\"ahler type and let $F$ be any fiber of $\pi$. Also, let $i_0$ and $i$ be the respective immersions of the fibres in $X$.  As $B$ is arcwise connected, from Ehresmann's theorem, it follows that $F$ and $F_0$ have the same homotopy type, which implies the exact sequence
$$
\begin{CD}
0 @>>> H^1(B) @ >\pi^*>> H^1(X) @>i^*>> H^1(F).
\end{CD}
$$
But Vaisman proved (\cite{va}) that if a compact LCK manifold of dimension at least $2$ is of K\"ahler type, then the LCK metric is actually GCK. Hence, if $F_0$ has dimension at least $2$,  it follows that $i_0^*([\theta])=0.$ From the exact sequence above, we see $[\theta]$ is a pull-back, and hence the above Lemma implies that  $g$ is GCK.

\medskip
To prove {\bf 2)}, observe that if the genus of $F$ is $0$ then $\pi^*$ is an isomorphism between $H^1(B)$ and $H^1(X)$, so the Lemma applies again. 

If the genus is at least $2$, we argue as follows. First, by Uniformization Theorem, after a  conformal change of $g$ we may assume $g_{\vert F}$ has (negative) constant curvature. On the other hand, by \cite{pv} we get that $[\theta]_{\vert F}$ is the Poincar\'e dual of the character $\chi$ of $g_{\vert F}$ (see \eqref{ch} for the definition of the character). But this character is trivial, since the Riemannian universal cover of $(F, g_{\vert F})$ is the Poincar\'e half-plane with the metric of negative constant curvature whith respect to which every homothety is an isometry. Hence $[\theta]_{\vert F}=0$, and so again $[\theta]$ is a pull-back from $B.$

\medskip

We now prove {\bf 3)}. As the $j-$invariant is a holomorphic map and $B$ is compact, we see that all the fibers are isomorphic and hence, by Fischer-Grauert theorem \cite{fg}, the map $\pi$ is a locally trivial fibration. Hence $X$ has a finite cover $X'$ which is a principal elliptic bundle  and the fiber $F$ acts holomorphically on $X'$. In particular, every $S^1\subset F$ acts holomorphically on $X'$. Hence, if $X'$ is not GCK, then by a result in  \cite{ov3} there exists a K\"ahler covering with K\"ahler metric given by an automorphic global potential. On the other hand, compact LCK manifolds with such coverings where shown to be complex deformations of Vaisman manifolds \cite{ov4}. \endproof

\hfill

Finally, let us give the

\hfill

\noindent{\bf Proof of Corollary 1.} Assume the product $X_1\times X_2$ has a LCK structure. As $X_1$ and $ X_2$ are regular compact Vaisman manifolds, they are total spaces of  holomorphic submersions $\pi_i:X_1\ra B_i, i=1, 2$ onto (compact Hodge) manifolds $B_1, B_2$ with fibers elliptic curves $F_1, F_2$. But then 
$$\pi:X_1\times X_2\ra B_1\times B_2, \qquad \pi(x_1, x_2)=(\pi_1(x_1), \pi_2(x_2))$$
 is a holomorphic submersion with typical fiber $F_1\times F_2$ which is a $2$-dimensional torus and is  of K\"ahler type. As the first Betti number of a compact Vaisman manifold is odd,  $b_1(X_1\times X_2) $ is even and we see that $X_1\times X_2$ is not biholomorphic to a Vaisman manifold. Then, from the above Theorem it follows that $X_1\times X_2$ is of K\"ahler type. But this forces $X_1, X_2$ to be of K\"ahler type as well, which is absurd.\endproof

\medskip

\section{Appendix on elliptic bundles and elliptic curves}
\noindent{\bf A.} Let us first recall some facts about elliptic bundles. Fix a genus 
one curve $E$ and let $E_0=(E, O)$ be the elliptic curve obtained by fixing 
some arbitrary point $O\in E.$ Fixing $O$ allow us to give a group 
structure on $E$. The group of automorphisms  $\Aut(E)$ is given by the 
extension
$$0\ra \mathrm{Trans}(E) \ra \Aut(E)\ra \Aut(E_0)\ra 0$$
where $\mathrm{Trans}(E)$ is the subgroup of $\Aut(E)$ given by translations 
and $\Aut(E_0)$ is the group of automorphisms of $E$ fixing $O.$ Now the 
group $\Aut(E_0)$ is usually $\Z_2$ (and consists of the antipodal map 
$x\mapsto -x$) except for some cases when $\Aut(E_0)$ is finite of 
order $4$ or $6$ (these particular kind of elliptic curves are 
called ``curves with complex multiplication''). See \cite[p. 143]{BPV}.

For an arbitrary complex manifold $M$ let $\mathrm{PBun}_E(M)$ respectively 
$\mathrm{Bun}_E(M)$ be the set of principal bundles, respectively the set of  
elliptic bundles on $M$ with fiber $E.$ The above exact sequence implies 
$$0\ra
\mathrm{PBun}_E(M)
\ra
\mathrm{Bun}_E(M)
\ra
H^1(M, \Aut(E_0))
$$
See again \cite[p. 143]{BPV}.
 
As $\Aut(E_0)$ is a finite group, we see that for any elliptic bundle 
$X\ra M$ there is a finite cover $M'$ of $M$ such that $X'=X\times_M M'$ 
is a principal bundle, in other words any elliptic bundle has a finite cover 
which is an  elliptic principal bundle (see \cite[p. 147]{BPV}).

\medskip

\noindent{\bf B.} We now recall some classical facts about the $j$-invariant. Let $\cae=\C/\langle 1,\tau\rangle$, $\tau\in \C$, $\Im \tau>0$ be a framed elliptic curve. Its $j$-invariant is the complex function
$$j(\cae)=j(\tau)=1728\frac{g_2^3(\cae)}{\Delta(\cae)},$$
where
\begin{equation*}
\begin{split}
g_2(\cae)&=g_2(\tau)=60\sum_{(m,n)\in\Z\setminus\{0\}}(m+n\tau)^{-4},\\
g_3(\cae)&=g_3(\tau)=160 \sum_{(m,n)\in\Z\setminus\{0\}}(m+n\tau)^{-6},
\end{split}
\end{equation*}
these two series being known to be absolutely convergent, and 
$$\Delta(\cae)=g_2^3(\cae)-27g_4^2(\cae).$$
From the very definition, the $j$-invariant is an analytic function of $\tau$.

If now $\pi:X\ra T$ is an analytic family of elliptic curves, one defines $J:T\ra\C$ by $J(t)=j(\cae_t)$. As $T$ is a manifold and hence locally simply connected, we may suppose the analytic family to be analytically framed. This implies that $J$ is analytic, as a composition of the analytic maps $\tau\mapsto j(\tau)$ and the period map $t\mapsto\tau(t)$.

{\small

\noindent {\sc Liviu Ornea\\
University of Bucharest, Faculty of Mathematics, \\14
Academiei str.,010014 Bucharest, Romania. \emph{and}\\
Institute of Mathematics ``Simion Stoilow" of the Romanian Academy,\\
21, Calea Grivitei Street
010702-Bucharest, Romania }\\
\tt liviu.ornea@imar.ro, \ \ lornea@gta.math.unibuc.ro

\hfill

\noindent {\sc Maurizio Parton\\
{\sc  Dipartimento di Scienze, Universita di Chieti--Pescara, \\
Pescara, Italy }\\
\tt parton@sci.unich.it

\hfill

\noindent {\sc Victor Vuletescu\\ University of Bucharest, Faculty of
Mathematics,
\\14
Academiei str., 010014 Bucharest, Romania.}\\
\tt vuli@gta.math.unibuc.ro
}
}

\end{document}